\documentclass{article}

\usepackage{amsmath,amsfonts,amssymb,amsthm}
 \newtheorem{theorem}{Theorem}[section]
 \newtheorem{corollary}[theorem]{Corollary}
 \newtheorem{lemma}[theorem]{Lemma}
 
 \theoremstyle{definition}
 \newtheorem{defn}[theorem]{Definition}
 \theoremstyle{remark}
 \newtheorem{remark}[theorem]{Remark}
 \numberwithin{equation}{section}

\newcommand{\Z}{\mathbb{Z}}

\newcommand{\M}{\mathbb{M}}

\newcommand{\CC}{\mathbb{C}}
\newcommand{\Q}{\mathbb{Q}}

\DeclareMathOperator{\Ker}{Ker}\DeclareMathOperator{\Irr}{Irr}\DeclareMathOperator{\Ch}{Ch}\DeclareMathOperator{\St}{St}\DeclareMathOperator{\Rad}{Rad}
\newcommand{\x}{\ensuremath{\textbf{x}}}\newcommand{\X}{\ensuremath{\textbf{X}}}

 \DeclareMathOperator{\Hom}{Hom}
\DeclareMathOperator{\diag}{diag}
\DeclareMathOperator{\La}{\mathrm{L}^{an}_D}
\begin{document}
\title{On the characters of pro-$p$ groups of finite rank}

\author{ A. Jaikin-Zapirain
\thanks{2000. \emph{Mathematics Subject Classification} Primary 20E18.
This work  has been partially supported by the MCYT Grants
BFM2001-0201, BFM2001-0180, FEDER and the Ram\'on y Cajal
Program.}\\
\footnotesize Departamento de Matem\'aticas,
\footnotesize Facultad de Ciencias,\\
\footnotesize Universidad Aut\'onoma de Madrid
\\}

\date{}
\maketitle

\section{Introduction}
We will say that a pro-$p$ group $G$ is {\bf perfect} if all
derived subgroups $G^{(k)}$ are open. Note that a pro-$p$ group is
perfect if and only if  $$\lambda_i(G)=|\{\lambda\in
\Irr(G)|\lambda(1)=p^i\}|$$ is finite for any $i\ge 0$. Here
$\Irr(G)$ is the set of the characters of the irreducible smooth
complex representations of $G$. When $G$ is uniform (a
torsion-free powerful pro-$p$ group), $G$ is perfect if and only
if  $[G,G]$ is open in $G$.

In this note we shall investigate  the function
$$\zeta^{ch}_G(s)=\sum_{i=0}^\infty \lambda_i(G)p^{-is}=\sum_{\lambda\in
Irr(G)}|\lambda(1)|^{-s},$$ when $G$ is perfect $p$-adic analytic
pro-$p$ group. The main result is the following theorem.
 \begin{theorem}\label{main}
Let $G$ be a perfect $p$-adic analytic pro-$p$ group. Then
$\zeta^{ch}_G(s)$ is a rational function in  $p^s$ if one of the
following conditions holds:

1. $G$ is a  uniform pro-$p$ group or

2. $p>2$.
\end{theorem}
The proof of this theorem is based on the correspondence between
the characters of a uniform pro-$p$ group and the orbits of the
action of the group on the dual of its Lie algebra. This result is
an analogue of the Kirillov theory, introduced first in the
context of nilpotent Lie groups and then used in many other
situations (see \cite{Ki}). The correspondence is quite explicit
and it also gives the exact formula for characters in some cases.
It permits us to  obtain a stronger version of Theorem \ref{main}
in these cases.

We introduce a generalization of the notion of a powerful pro
$p$-group: we say that a pro $p$-group $G$ is {\bf $k$-powerful}
if $[G,G]\le G^{p^k}$. We say that a pro-$p$ group is {\bf
$k$-uniform} if $G$ is $k$-powerful without torsion. We have the
following generalization of Theorem \ref{main}.
\begin{theorem}\label{main2} Let $N$ be
 a 1-uniform (uniform) pro-$p$ group if $p\ge 5$  and $2$-uniform pro-$p$ group
if $p=3$. Then for any $g\in N$, $$\sum_{\lambda\in
Irr(G)}\lambda(g)|\lambda(1)|^{-s}$$ is a rational function in
$p^s$.
\end{theorem}
In Section \ref{howe} we describe Howe's version of Kirillov's
correspondence. The correspondence permits us ``linearize'' the
problem and apply the potent tools of $p$-adic integration which
we introduce in  Section \ref{integr}. In Section \ref{uniform} we
prove Theorem \ref{main}(1) and Theorem \ref{main2}. Section
\ref{triples} is dedicated to some results from the character
theory of finite groups. We apply these results in Section
\ref{final}, where we prove Theorem \ref{main}(2).

 The notation is standard. If $G$ is a
profinite group, then $\Irr(G)$ denotes the set of the (complex)
irreducible smooth (not only linear) characters of $G$. If $L$ is
a $\Z_p$-Lie algebra of finite rank, then $\Irr(L)$ is the set of
the irreducible characters of the additive group of $L$ and
$L^*=\Hom_{\Z_p}(L,\Z_p)$ is the dual of $L$. The set $\Ch(G)$ is
the set of the admissible smooth characters of $G$. So any element
of $\Ch(G)$ is a finite sum of elements from $\Irr(G)$. The scalar
product of character we denote by $<\ ,\ >$. If $N$ is a normal
group of finite index and $\chi\in \Irr(N)$, then
$$\Irr(G|\chi)=\{\lambda\in \Irr(G)|<\lambda_N,\chi>\ne 0\}$$ and
$\Ch(G|\chi)$ is the set of the finite sums of elements from
$\Irr(G|\chi)$. We will say that characters from $\Ch(G|\chi)$ lie
over $\chi$. We will use  $[\ ,\ ]_L$ to denote the Lie brackets
and $[\ ,\ ]_G$ for the group commutator. If $g$ is an elemethe of
a group, $o(g)$ will mean the order of element $g$. The order of a
non-zero complex number is its order in $\CC^*$.

\section{Correspondence between characters and coadjoint
orbits}\label{howe} In this section we describe the correspondence
between characters of a uniform pro-$p$ group $N$ and coadjoint
orbits of action of $N$ on the dual of the Lie algebra associated
with $N$. This section is based on the paper \cite{Ho}, where
almost all proofs of the results of this section can be found (at
least in the case $p>2$).

Remind that if $N$ is a uniform pro-$p$ group, then from Lazard
(see, for example, \cite{DDMS}) we know that we can associate a
$\Z_p$-Lie algebra $L$ with $N$. This Lie algebra can  be
identified with $N$ as a set and the Lie operations are defined by
$$g+h=\lim_{n\to \infty} (g^{p^n}h^{p^n})^{p^{-n}}, \ [g,h]_L=\lim_{n\to \infty}[ g^{p^n}, h^{p^n }]_G^{p^{-2n}}.$$
Since $(L,+)$  is a pro-$p$ group, we can consider the set
$\Irr(L)$. Note that $N$ acts on $L$ by conjugation, whence $N$
also acts on $\Irr(L)$. It is easy to see that all $N$-orbits of
the last action are finite.

\begin{theorem}\label{charunif} Let $N$ be
 a 1-uniform (uniform) pro-$p$ group if $p\ge 5$  and $2$-uniform pro-$p$ group
if $p=3$ and let $\Omega$ be a $N$-orbit in $\Irr(L)$. Define
$\Phi_{\Omega}\colon N\to \CC$ by means of
$$\Phi_\Omega(u)=|\Omega|^{-\frac{1}{2}}\sum_{\omega\in \Omega}\omega(u).$$
Then $\Phi_\Omega\in \Irr(N)$ and all characters of $G$ have this
form. Moreover, two different orbits $\Omega_1$ and $\Omega_2$
give two different $\Phi_{\Omega_1}$ and $\Phi_{\Omega_2}$.
\end{theorem}
 If $N$ satisfies the
hypothesis of the previous theorem and it is also perfect, then we
obtain that $\lambda_i(N)$ is equal to the number of $N$-orbits in
$\Irr(L)$ of size $p^{2i}$.  Next theorem shows that this is true
for all uniform pro-$p$ groups.
\begin{theorem}\label{charunif2} Let $N$ be  a uniform pro-$p$
group, $\alpha\in \Irr(N)$ and $\Omega$  a $N$-orbit in $\Irr(L)$.
If $<\alpha, \Phi_\Omega>\ne 0$, then $\alpha(1)=\Phi_\Omega(1)$.
In particular, if $N$ is perfect, then $\lambda_i(G)$ is equal to
the number of $N$-orbits in $\Irr(L)$ of size $p^{2i}$.
\end{theorem}

If $w\in \Irr(L,+)$, we define
$$B_w(l,k)=w([l,k]_L),\ l,k\in L.$$
Then  $B_w$ is a bilinear form on $L$. Put $\Rad(w)=\{l\in
L|B_w(l,L)=1\}$. An important point is that $\Rad(w)$ is equal to
$\St_N(w)=\{g\in N|w^g=w\}$. Hence by the previous theorem, we
have that if $N$ is a uniform pro-$p$ group then
$$\zeta^{ch}_N(s)=\sum_{w\in Irr(L)}|L:\Rad(w)|^{\frac{1}{2}(-s-2)}.$$ Therefore we
are interested in  the function $\zeta(s)= \sum_{w\in
Irr(L)}|L:\Rad(w)|^{-s}.$

\section{P-adic integration}\label{integr}
In this section we give  an explanation of the notion of $p$-adic
integral and we introduce the facts that we will use later. More
detailed discussion of this subject can be found in
\cite{De,DvD,duS}.

We use the standard notation for $p$-adic sets. $|\ |$ is the
standard $p$-adic valuation on $\Q_p$: if $a \in p^k\Z_p\setminus
p^{k+1}\Z_p$, entonces $|a|=p^{-k}$. $\mu$ will be the Haar
mesure on $\Q_p^n$. We always suppose that $\mu(\Z_p^n)=1$.

Let $\X=(X_1,\cdots,X_M)$ be  $M$ commuting indeterminates and let
$\Q_p[[\X] ]$ denote the set of formal power series over $\Q_p$.
We define the following subsets of $\Q_p[[\X] ]$:

1. $\Z_p[[\X] ]$ denotes the set of power series over $\Z_p$;

2. $\Q_p\{\X\}$ consists of all formal  power series
$\sum_ia_i\X^i$ such that $|a_i|\to 0$ as $|i|\to \infty$.

3. $\Z_p\{\X\} =\Z_p[[\X] ]\cap \Q_p\{\X\}.$

We define the function $D\colon \Z_p^2\to \Z_p$ by
$$D(x,y)=\left \{ \begin{array}{ll}
x/y & \textrm{\ if \ } |x|\le |y| \textrm{\ and\ }y\ne 0, \\
0 & \textrm{\ otherwise.\ }
\end{array}\right .$$
For $n>0$ we define $P_n$ to be the set of nonzero $n$th powers in $\Z_p$.

We define the following language considered in \cite{DvD}:

Let $\La$ be the language with logical symbol
$\forall,\vee,\wedge,=$, a contable number of variables $X_i$ and
\begin{enumerate}
\item an $m$-place operation symbol $F$ for each $F(\X)\in\Z_p\{\X\}$, $m\ge 0$;
\item a binary operation symbol $D$;
\item a unitary relation symbol $P_n$ for each $n\ge 2$.
\end{enumerate}

Each formula $\phi(x_1,\cdots,x_M)$ in the language $\La$ defines a subset
$$M_\phi=\{\x\in \Z_p^M|\phi(\x) \textrm{\ is true in\ } \Z_p\},$$
where we interpret
\begin{enumerate}
\item each $F\in\Z_p\{\X\}$, as a function $f\colon \Z_p^n\to \Z_p$ defined by $f(\x)=F(\x)$;
\item the binary operation symbol $D$ as the function $D$;
\item $P_n(x)$ to be true if $x\in P_n$.
\end{enumerate}
 We call such subset $M_\phi$ definable.
 A function $f\colon V\to \Z_p$ is called definable if its graph is definable
 subset. Note that, in particular, a definable function is
 bounded.

With each pair of definable functions $f_1\colon \Z_p^M\to \Z_p$,
$f_2\colon \Z_p^M\to \Z_p$  and the definable subset $U$  of
$\Z_p^M$ we associate the following function:
$$I(f_1,f_2,U,s)=\int_U|f_1(\x)|^s|f_2(\x)|d\mu.$$
We shall call this function definable integral.
\begin{theorem}\label{padicint}
Suppose that $I(f_1,f_2,U,s)$ is a definable integral. Then
\begin{enumerate}
\item $U$ is measurable,  and
\item $I(f_1,f_2,U,s)$ is a rational function in $p^{-s}$.
\end{enumerate}
\end{theorem}
\section{The uniform case}\label{uniform}
We begin this section with some known facts about endomorphisms of
$\Z_p^n$. Let $A, B\in \M_n(\Z_p)$ be  two matrices over $\Z_p$.
We write  $A\sim B$ if there are two invertible matrices $C_1$ and
$C_2$ over $\Z_p$ such that $C_1AC_2=B$.

If $A=(a_{ij})_{1\le i,j\le n}$ and $U\subseteq \{1,\cdots ,n\}$
we put $g_U(A)=\det [(a_{ij})_{i,j\in U}]$.

For any $1\le i\le n$ we fix
 an order on subsets of
$\{1,\cdots,n\}$ with $i$ elements. We put $h_0(A)=1$ and let for
$1\le i\le n$ $h_i(A)$ be $g_U(A)$ such that $|U|=i$ and for any
$U^\prime\subseteq \{1,\cdots,n\}$ with $|U^\prime|=i$ we have
$|g_{U^\prime}(A)|< |g_U(A)|$ or $|g_{U^\prime}(A)|= |g_U(A)|$
and $U\le U'$.

\begin{lemma} Let $A\in \M_n(\Z_p)$ be a matrix over $\Z_p$. Then

1. There are $s_1,\cdots,s_n\in \Z_p$, satisfying $|s_i|\ge
|s_{i+1}|$ for any $1\le i\le n$, such that
$A\sim\diag(s_1,\cdots,s_n)$. Moreover, if $A\sim\diag(t_1,\cdots,
t_n)$ with $|t_i|\ge |t_{i+1}|$ for any $1\le i\le n$, then
$|t_i|=|s_i|$.

2. For any $0\le i\le n-1$, $|D(h_i(A),h_{i-1}(A))|\ge
|D(h_{i+1}(A),h_i(A))|$ and
$$A\sim \diag
(h_1(A),D(h_2(A),h_1(A)),\cdots,D(h_n(A),h_{n-1}(A))).$$
\end{lemma}

\begin{proof}
It is enough to observe that $A\sim\diag(s_1,\cdots,s_n)$ if and
only if $\Z_p^n/A(\Z_p^n)\cong \oplus _{i=1}^n \Z_p/s_i\Z_p$.
\end{proof}
  It is not difficult to see that if $A\sim
B$ then $|h_i(A)|=|h_i(B)|$. Therefore, sometimes we will speak
about $|h_i(\phi)|$, where $\phi\in \Hom_{\Z_p}(M_1,M_2)$ and
$M_1\cong M_2\cong \Z_p^n$. Hence the previous lemma implies the
following result
\begin{lemma}
Let  $M_1\cong M_2\cong \Z_p^n$ be two $\Z_p$-modules, $0\ne z\in
\Z_p$ and $\phi\in \Hom_{\Z_p}(M_1,M_2)$. Then
$|M_1/\phi^{-1}(zM_2)|=$
$$ \left \{
\begin{array}{lll}
1 & \textrm{\ if \ } & |z|> |h_1(\phi)|,\\
|z|^{-k}|h_k(\phi)| &  \textrm{\ if } &
|D(h_{k}(\phi),h_{k-1}(\phi))|\ge |z|>
 |D(h_{k+1}(\phi),h_k(\phi))|,1\le k\le n-1\\
|z|^{-n}|h_n(\phi)| &  \textrm{\ if } &
|z|\le|D(h_{n}(\phi),h_{n-1}(\phi))|.
\end{array}\right.$$
\end{lemma}
\begin{proof} By the previous lemma, we can choose bases  of $M_1$  and $M_2$ such
that the matrix $A$ associated with $\phi$ in these bases is
diagonal and equal to  $$\diag
(h_1(A),D(h_2(A),h_1(A)),\cdots,D(h_n(A),h_{n-1}(A))).$$ Since
$|h_i(\phi)|=|h_i(A)|$ for all $i$, we obtain the lemma.
\end{proof}
Now, let $N$ be a uniform pro-$p$ group  and $L$ the Lie algebra
associated with $N$.  Put $L^*=\Hom_{\Z_p}(L,\Z_p)$. Since $N$
acts on $L$, $L^*$ has a structure of $N$-module.

Note that $\Irr(L)$ is the direct limit of $\Irr(L/p^kL)$.  For
each $i\ge 0$ we fix $\theta_i$ a $p^i$th primitive root of 1 in
$\CC$ and we construct a $N$-homomorphism $\Phi_i\colon L^*\to
\Irr(L/p^iL)$ in the following way:
$$\Phi_i(m)(l+p^iL)=\theta_i^{m(l)}.$$
The kernel of $\Phi_i$ is equal to $p^iL^*$ and so $\Phi_i$ is
surjective.
\begin{remark}
Note also that $L^*$  is  $N$-isomorphic to the inverse limit of
$(\Irr(L/p^kL),\phi_{k,i})$ where the homomorphism
$\phi_{k,i}\colon \Irr(L/p^kL)\to \Irr(L/p^iL)$ ($k\ge i$) is
defined by means of
$$\phi_{k,i}(m)(l+p^iL)=m({p^{k-i}}l+p^kL).$$
Therefore, $\{ \Phi_i\}$ can be also  constructed in a canonical
way.
\end{remark}

For any $0\ne z\in \Z_p$ such that $|z|=p^{-i}$ we define
$\Phi_z=\Phi_i$.

Let $\Psi\colon L^*\to \Hom_{\Z_p}(L,L^*)$ be the map defined in
the following way: if $l,k\in L$ and $m\in L^*$ then
$$\Psi(m)(l)(k)=(m([k,l]).$$
\begin{lemma} Let $m\in L^*$ and $0\ne z\in \Z_p$. Then
$$\Rad (\Phi_z(m))=(\Psi(m))^{-1}(zL^*).$$
\end{lemma}
\begin{proof} Let $|z|=p^{-i}$. Then we have the following
series of equivalent propositions:
\begin{eqnarray*}
 x\in \Rad(\Phi_z(m))\Leftrightarrow \Phi_z(m)([y,x])=1 \forall y\in L
\Leftrightarrow \\
\theta_i^{m([y,x])}=1  \forall y\in L\Leftrightarrow m([x,y])\in
z\Z_p  \forall y\in L \Leftrightarrow \\
\Psi(m)(x)(y)\in z\Z_p  \forall y\in L \Leftrightarrow
\Psi(m)(x)\in zL^* \Leftrightarrow x\in
(\Psi(m))^{-1}(zL^*).\end{eqnarray*}
\end{proof}
Let $\{e_1,\cdots,e_n\}$ be a basis of $L$ and
$\{f_1,\cdots,f_n\}$  a basis of $L^*$. So, any element $a$ from
$L$ or $L^*$ is identified with a vector $(a_1,\cdots,a_n)$ and we
can view $\Hom_{\Z_p}(L,L^*)$ as $\M_n(\Z_p)$. Then the entries of
$\Psi(a),a\in L^*$ are linear functions on $\{a_i\}$ and
$g_U(\Psi(a))$ are polynomials on $\{a_i\}$ for every $U\subseteq
\{1,\cdots,n\}$. This implies that $h_i(\Psi(a))$ are definable
functions in $\La$.

Define the sets $W= (L^*\setminus  pL^*)\times (p\Z_p\setminus
\{0\})$,
$$W_0=\{(a,z)\in W||z|> |h_1(\Psi(a))|\},$$
for any $1\le k\le n-1$,
 $$W_k=\{(a,z)\in W||D(h_{k}(\Psi(a)),h_{k-1}(\Psi(a)))|\ge |z|>
|D(h_{k+1}(\Psi(a)),h_k(\Psi(a)))|\}$$ and $$W_n=\{(a,z)\in
W||z|\le|D(h_{n}(\Psi(a)),h_{n-1}(\Psi(a)))|\}.$$ If $(a,z)\in
W_k$ define  $\alpha(a,z)=D(z^k,h_k(\Psi(a)))$. Then $\alpha$ is
a definable function on $W$.
 Using two previous lemmas, we obtain that
\begin{corollary} Let $(a,z)\in W$. Then
$|L:Rad(\Phi_z(a))|=|\alpha(a,z)|^{-1}$.
\end{corollary}

Now, suppose that $N$ is perfect. In this case we have the next
important property:
 \begin{lemma} \label{perfect}
 Let $N$ be  a perfect uniform
pro-$p$ group and $L$ the Lie algebra associated  with $N$. Then
the number
$$m(L)=\min\{k|p^kL\subseteq [L,L]_L\}$$
is finite and for any $w\in Irr(L)$ we have $o(w)\le
|L:\Rad(w)|p^{m(L)}$.
\end{lemma}
\begin{proof} The finiteness of $m(L)$ follows from the
perfectness of $L$.

Since $p^{m(L)}L\le [L,L]_L$, $|L:\Rad(w)|p^{m(L)}L\le
[L,\Rad(w)]_L$. Hence $$w(|L:\Rad(w)|p^{m(L)}L)=1,$$ and so
$|L:\Rad(w)|p^{m(L)}\ge o(w)$.
\end{proof}
\begin{corollary}\label{corper}
Let $(a,z)\in W$. Then
$$ |p^{m(L)}\alpha(a,z)|\le |z|.$$
\end{corollary}
\begin{proof}Note that if $w=\Phi_z(a)$, then $|z|=o(w)^{-1}$.
\end{proof}
\begin{theorem}\label{unif}
Let $N$ be a perfect uniform pro-$p$ group. Then $\zeta^{ch}_N(s)$
is a rational function in $p^s$.
\end{theorem}
\begin{proof} Let $L$ be the Lie algebra associated with $N$ and $n$ the $\Z_p$-rank of $L$.
First note that if $(a,z)\in W$ and $w=\Phi_z(a)$, then
$|z|=o(w)^{-1}$. Hence  $$\mu(\{(a,z)\in
W|\Phi_z(a)=w\}=(p-1)p^{-1}o(w)^{-(n+1)}.$$ Therefore
\begin{eqnarray*}\zeta(s)-1=\sum_{1_L\ne w\in Irr(L)}|L:Rad
(w)|^{-s}=\\ \sum_{1_L\ne w\in \Irr(L)}\int_{ (a,t)\in
W,\Phi_z(a)=w} \frac{p(p-1)^{-1} }{|z|^{(n+1)}}|L:Rad(\Phi_z(a)|^{-s}dadz=\\
p(p-1)^{-1}\int_{W}\frac{1}{|z|^{(n+1)}}|L:Rad(\Phi_z(a)|^{-s}dadz=\\
p(p-1)^{-1}\int_{W}
 |z|^{-(n+1)}|\alpha(a,z)|^{s}dadz.
\end{eqnarray*}
By Corollary \ref{corper} we have that if  $(a,z)\in W$, then
\begin{eqnarray*}|z|^{-1}=|z|^{-1}|p^{m(L)}\alpha(a,z)||p^{m(L)}\alpha(a,z)|^{-1}
=\\
|D(p^{m(L)}\alpha(a,z),z)||p^{m(L)}\alpha(a,z)|^{-1}.\end{eqnarray*}
Therefore we have:
\begin{eqnarray*}
\int_{W}
 |z|^{-(n+1)}|\alpha(a,z)|^{s}dadz=\\ p^{m(L)(n+1)}(\int_{W}
 |D(p^{m(L)}\alpha(a,z),z)|^{n+1}|\alpha(a,z)|^{s-n-1}dadz
 )
\end{eqnarray*}

Now from Theorem \ref{padicint} we obtain that the last integral
is a rational function in $p^s$, whence $\zeta(s)$ is a rational
function in $p^s$. Note that since $\zeta(s)=\sum a_ip^{-2is}$, we
have that $\zeta(s)$ is a rational function in  $p^{2s}$. We
conclude from the last paragraph of Section 2 that
$\zeta^{ch}_N(s)$ is a rational function in $p^s$.
\end{proof}
Let $N$ be a perfect uniform pro-$p$ group, $g\in N$ and let
$\mu_i$ be the sum of all irreducible characters of $N$ of degree
$p^i$. Then
$$\sum_{\lambda\in
Irr(G)}\lambda(g)|\lambda(1)|^{-s}=\sum_{i=0}^{\infty}\mu_i(g)p^{-si}.$$
Therefore the following theorem implies Theorem \ref{main2}.
\begin{theorem} Let $N$ be
 a 1-uniform (uniform) pro-$p$ group if $p\ge 5$  and $2$-uniform pro-$p$ group
if $p=3$. Then for any $g\in N$,
$$\sum_{i=0}^{\infty}\mu_i(g)p^{-si}$$ is a rational function in
$p^s$.
\end{theorem}
\begin{proof}
By Theorem \ref{charunif},
$$\mu_i=\sum_{|\Omega|=p^{2i}}\Phi_\Omega=p^{-i}\sum_{w\in
\Irr(L),\ |L:\Rad(w)|=p^{2i}}w.$$ Hence
$$\sum_{i=0}^{\infty}\mu_i(g)p^{-si}=\sum_{w\in \Irr(L)}w(g)|L:\Rad(w)|^{-(s+2)/2}.$$
 Let $\theta_m$ be a $p^m$th
primitive root of 1. If the order of $w$ is equal to $p^m\ge p$,
then
$$\sum_{\sigma\in \textrm{Gal}
(\Q(\theta_m)/\Q)}w^\sigma(g)=\left\{\begin{array}{lll}
(p-1)p^{m-1} &
\textrm{\ if \ } & g\in \Ker w\\
-p^{m-1} & \textrm{\ if \ } & o(w(g))=p \\
0 & \textrm{\ if \ } & o(w(g))>p .
\end{array}\right.$$
Therefore we have
$$\begin{array}{lll}
\sum_{i=0}^{\infty}\mu_i(g)p^{-si} & = & \sum_{w\in \Irr(L),\
o(w(g))=1}|L:\Rad(w)|^{-(s+2)/2}-\\ & & \frac{1}{p-1}\sum_{w\in
\Irr(L),\ o(w(g))=p}|L:\Rad(w)|^{-(s+2)/2}.
\end{array}$$ Now, note that the sets  $W_1=\{(a,z)\in
W|a(g)\equiv 0 (\mod z)\}$ and $W_2=\{(a,z)\in W|pa(g)\equiv 0
(\mod z)\}$ are definable. Following the proof of the previous
theorem, we obtain that
$$\sum_{w\in \Irr(L),\
o(w(g))=1}|L:\Rad(w)|^{-s}-1=p(p-1)^{-1}\int_{W_1}|z|^{-(n+1)}|\alpha(a,z)|^{s}dadz
,$$ and
 $$\sum_{w\in
\Irr(L), \
o(w(g))=p}|L:\Rad(w)|^{-s}=p(p-1)^{-1}\int_{W_2\setminus W_1}
 |z|^{-(n+1)}|\alpha(a,z)|^{s}dadz.$$
 Hence, as in the previous theorem, we conclude that
 $\sum_{i=0}^{\infty}\mu_i(g)p^{-si}$ is a rational function in
 $p^s$.
\end{proof}

\section{Character triples}\label{triples}
Let $G$ be a group, $N$ a normal subgroup of $G$ of finite index
and $\theta \in Irr(N)$ $G$-invariant irreducible character of
$N$. Under these hypotheses we say that $(G,N,\theta)$ is a {\bf
character triple}.  We refer the reader to \cite[Section 11]{Is}
for general discussion on character triples. The main idea consist
in replacing  $(G,N,\theta)$ by another character triple $(\Gamma,
A, \lambda)$ in which $\Gamma/A\cong G/N$ and $\lambda$ is linear.
Moreover, the character theory of extensions of $\lambda$ on
$\Gamma$ and $\theta$ on $G$ is the same. (The next definition
explains what this means exactly.) For instance, we will have that
$\lambda^\Gamma$ and $\theta^G$ have the same numbers of
irreducible constituents with the same ramifications.
\begin{defn}\label{isomtriple}
Let $(G,N,\theta)$ and $(H,M,\phi) $ be character triples and
$\tau\colon H/M\to G/N$ be an isomorphism. For $M\le T\le H$, let
$T^\tau$ denote the inverse image in $G$ of $\tau(T/M)$. For every
such $T$, suppose there exists a map $\delta_T\colon
\Ch(T|\theta)\to \Ch(T^\tau|\phi)$ such that the following
conditions hold for $T$, $K$ with $M\le K\le T\le H$ and
$\chi,\psi\in \Ch(T|\theta)$.
\begin{enumerate}
\item
$\delta_T(\chi+\psi)=\delta_T(\chi)+\delta_T(\psi)$;
\item
$\langle \chi,\psi\rangle=\langle
\delta_T(\chi),\delta_T(\psi)\rangle$;
\item
$\delta_K(\chi_K)=(\delta_T(\chi)_{K^\tau})$;
\item
$\delta_T(\chi\beta)=\delta_T(\chi)\beta^\tau$ for $\beta\in
Irr(T/M)$.
\end{enumerate}
Let $\delta$ denote the union of the maps $\delta_T$. Then $(\tau,
\delta)$ is an {\bf isomorphism} from $(G,N,\theta)$ to
$(H,M,\phi)$.
\end{defn}

Now, let $p>2$ and $G$ be a $p$-adic analytic pro-$p$ group. We
know that $G$ has a $2$-uniform subgroup $N$ of finite index. Let
$L$ be the Lie algebra associated with $N$.  Let $\chi$ be an
irreducible character of $N$. Then, from Theorem \ref{charunif},
it follows that there exists
  a $G$-orbit $\Omega$  in $\Irr(L)$, such that
$\chi=\Phi_{\Omega}$. Let $w\in \Omega$. By \cite[Theorem
6.11]{Is}, in order to understand the character theory of
extensions of  $\chi$ on $G$ we should investigate the character
theory of extensions of  $\chi$ on $I=I_G(\chi)=\{g\in
G|\chi^g=\chi\}$. This problem is solved in the next theorem.
\begin{theorem}\label{triplesreduction}
Let $G$, $N$, $I$, $L$, $\chi$, $\Omega$ and $w$ as before. Then
$w_{\St_N(w)}\in \Irr(\St_N(w))$ and  $(I, N,\chi)$ is isomorphic
to $(\St_G(w),\St_N(w),w_{\St_N(w)})$.
\end{theorem}
This result is practically Proposition 1.2 from \cite{Ho}. But,
since we use a different notation, we present a sketch of the
proof.
\begin{proof}
First, by \cite[Lemma 1.4]{Ho}, there exists a subgroup $A$ of $N$
such that
\begin{enumerate}
\item
A is a uniform pro-$p$ group and so it is also a Lie subalgebra of
$L$;
 \item  $B_w(A,A)=1$ and $A$ is maximal with this property;
\item $w_A$ is also a groups character of $A$ respect
multiplication and $\chi=(w_A)^N$; \item $A$ is
$\St_G(w)$-invariant.
\end{enumerate}
Put $S=\St_G(w)A$. Then we have that $I=SN$ and $A=S\cap N$. Note
that $w_A$ is $S$-invariant. First, we will see that $(I, N,\chi)$
and  $(S,A,w_A)$ are isomorphic.

Let $\tau$ be the natural isomorphism between $S/A=S/(S\cap N)$
and $I/N=SN/N$, $T$ a subgroup of $S$, which contains $A$, and
$T^\tau$ the inverse image in $I$ of $\tau(T/A)$. Hence
$T^\tau=TN$. Define $\delta_T\colon \Ch(T|w_A)\to
\Ch(T^\tau|\chi)$ by means of
$$\delta_T(\alpha)=\alpha^{T^\tau},\ \alpha\in \Ch(T|w_A).$$
Since $w_A^N=\chi$ and $N$ fixes $\chi$, $\delta_T$ sends $
Ch(T|w_A)$ to $Ch(T^\tau|\chi)$.

We have to check four properties of Definition \ref{isomtriple}.
The fist one is clear.
 In order to prove the second one, observe that
 $$(w_A)^{T^\tau}=((w_A)^N)^{T^\tau}=\chi^{T^\tau}.$$
 In particular,
 $$<(w_A)^{T^\tau},(w_A)^{T^\tau}>=<\chi^{T^\tau},\chi^{T^\tau}>=<\chi,(\chi^{T^\tau})_N>=|T^\tau/N|.$$
 On the other hand,
 $$(w_A)^{T^\tau}=((w_A)^T)^{T^\tau}=( \sum_{\lambda\in
 \Irr(T|w_A)}{\lambda(1)}\lambda)^{T^\tau}=\sum_{\lambda\in
 \Irr(T|w_A)}{\lambda(1)}\lambda^{T^\tau}.$$
 Hence we have
 $$\begin{array}{lll}
 |T^\tau/N| &=
 &<(w_A)^{T^\tau},(w_A)^{T^\tau}>=\sum_{\lambda,\mu\in\Irr(T|w_A)}{\lambda(1)}\mu(1)<\lambda^{T^\tau},\mu^{T^\tau}>\le\\
 & &\sum_{\lambda\in\Irr(T|w_A)}{\lambda(1)}^2=|T/A|=|T^\tau/N|.
 \end{array}$$
This implies that
$$<\lambda^{T^\tau},\mu^{T^\tau}>=\left \{ \begin{array}{lll}
1 & \textrm{if} & \lambda=\mu\\
0 & \textrm{if} & \lambda \ne \mu \end{array}\right.$$
 But this is exactly the second condition.

The third property is a direct consequence of \cite[Problem
5.2]{Is} and the forth one can be obtained from \cite[Problem
5.3]{Is}, bearing in mind that $\beta^\tau_T=\beta$.

Now we will see that $(S,A,w_A)$ and  $(\St_G(w),\St_N(w),
w_{\St_N(w)})$ are isomorphic. Since $\St_N(w)=\St_G(w)\cap A$ we
have that $S/A$ and $\St_G(w)/\St_N(w)$ are isomorphic. Let $\tau$
be the natural isomorphism between $S/A$ and $\St_G(w)/\St_N(w)$,
$T$ a subgroup of $S$, which contains $A$, and $T^\tau$ the
inverse image in $I$ of $\tau(T/A)$. Hence $T^\tau=T\cap
\St_G(w)$. Define $\delta_T\colon \Ch(T|w_A)\to
\Ch(T^\tau|w_{\St_G(w)})$ by means of
$$\delta_T(\alpha)=\alpha_{T^\tau},\ \alpha\in \Ch(T|w_A).$$
We  left to the reader to check that this map is an isomorphism of
triples.
 \end{proof}
Let $Q$ be a finite group and $F$ be a free  group on $|Q|$
variables (so $F$ is generated by $x_q$, $q\in Q$).  Define an
homomorphism $\phi\colon F\to Q$ by means of $\phi(x_q)=q$. Let
 $H$ be the kernel of this homomorphism.  Put
$\bar F=F/[H,F]$ and $\bar H=H/[H,F]$. Then $\bar H$ is an abelian
group of finite rank. Hence we can write $\bar H=A\oplus B$ where
$A$ is torsion-free and $B$ is a torsion group. Let $\beta$ and
$\alpha$ be $F$-invariant linear characters of $H$. Since they are
$F$-invariant, we can see $\alpha$ and $\beta$ as characters of
$\bar H$. We will need the following criterion.
\begin{lemma} With the previous notation suppose that for every $b\in B\cap [\bar F,\bar F]$ the orders of $\alpha(b)$ and
$\beta(b)$ are same. Then the triples $(F, H, \alpha)$ and $(F,
H,\beta)$ are isomorphic.
\end{lemma}
\begin{proof} We split the proof in a number of
steps.

 {\it Step 1. }
 Let $\gamma$ be a linear character of $\bar H$ which is trivial on
 $B\cap [\bar F,\bar F]$. Then there a character $\tau$ of $F$ such that
 $\tau_H=\gamma$.

 Note that the comutator of $\bar F$ is finite,
 so $[\bar F,\bar F]\cap \bar H= B\cap [\bar F,\bar F]$. Hence
 $\gamma $ is also trivial on $[\bar F,\bar F]\cap \bar H$.
 Therefore we can see $\gamma$ as a character of $C=\bar H/([\bar F,\bar F]\cap \bar
 H)$. But $C$ can be seen as a subgroup of $\bar F/[\bar F,\bar
 F]$ and so we can extend $\gamma$ on $\bar F/[\bar F,\bar
 F]$.

{\it Step 2. } Let $\gamma$ be a linear character of $\bar H$
which is trivial on
 $B\cap [\bar F,\bar F]$. Then $(F, H, \alpha)$ and $(F, H, \alpha\gamma)$ are
 isomorphic.

 We will indicate the isomorphism. By the previous step there exists $\tau\in \Irr(F)$ which extends $\gamma$.
 Let $T$ be a subgroup of $F$
 which contains $H$ and $\chi$ a character of $T$ lying over
 $\alpha$.  We put $\delta_T(\chi)=\chi\tau_T$. It is not
 difficult to see that it is really isomorphism of triples.

{\it Step 3. } Final step.

The conditions of lemma imply that $\alpha_{B\cap [\bar F,\bar
F]}$ and $\beta_{B\cap [\bar F,\bar F]}$ are conjugated by an
element $\sigma$ of the Galois group of $\Q$. Hence
$\alpha=\beta^{\sigma}\gamma$, where $\gamma$ is a linear
character of $\bar H$ which is trivial on
 $B\cap [\bar F,\bar F]$.

By the previous step, $(F, H, \alpha)$ and $(F, H,
\beta^{\sigma})$ are
 isomorphic. On the other hand, $(F, H, \beta)$ and $(F, H,
\beta^{\sigma})$ are also
 isomorphic.
\end{proof}

Now let $\{s_i=s_i(x_q|q\in Q)\}\subset [F,F],$ be a finite set of
elements of $F$ such that $B\cap [\bar F,\bar F]=\{s_i[H,F]\}$. We
denote by $S_Q$ the set $\{s_i\}$.

Now, let $N\to G \overset \phi  \to Q$ be an extension of $Q$,
i.e. $\phi$ is a surjective homomorphism with kernel $N$. Let
$\{t_q|q\in \Q\}$ be a transversal for $N$ in $G$, such that
$\phi(t_q)=q$. Let $\alpha$ be a $G$-invariant linear character of
$N$. Define a vector
$$S_Q(G,N,\phi,
\alpha)=(o(\alpha(s_1(t_q|q\in \Q))),\ldots,o(\alpha( s_k(t_q|q\in
\Q)))),$$
 where $k=|S_Q|$ and $o(r)$ means the order of $r$.
\begin{corollary}\label{crit}
\begin{enumerate}
\item The vector $S_Q(G,N,\phi,
\alpha)$ does not depend on the choice of transversal for $N$ in
$G$.
\item
There exists only finite number of posilibilities for
$S_Q(G,N,\phi, \alpha)$.
\item If $N_1\to G_1 \overset {\phi_1}  \to Q_1$ and $N_2\to G_2 \overset {\phi_2}  \to Q_2$ are two  extensions of
$Q$, $\alpha_i$ is a $G_i$-invariant linear character of $N_i$ for
$i=1,2$ and $S_Q(G_1,N_1,\phi_1, \alpha_1)=S_Q(G_2,N_2,\phi_2,
\alpha_2)$, then $(G_1,N_1, \alpha_1)$ and $(G_2,N_2, \alpha_2)$
are isomorphic.
\end{enumerate}
\end{corollary}

\section{The general case}\label{final}
Throughout of this section we suppose that $p>2$ and $G$ be a
$p$-adic analytic pro-$p$. We can find an open normal 2-uniform
subgroup $N$ of $G$.  Let $L$ be the Lie algebra associated with
$N$. Since $G$ acts on $L$, $G$ also acts on $L^*$. We will use
the notation of Section \ref{uniform}.  Let $\{e_1,\cdots,e_n\}$
be a basis of $L$ and $\{f_1,\cdots,f_n\}$  a basis of $L^*$. So,
any element $a$ from $L$ or $L^*$ is identified with a vector
$(a_1,\cdots,a_n)$ from $\Z_p^n$.

For any subgroup $N\le K\le G$, define the following set
 $$\Irr(L)_K=\{w\in \Irr(L)|N\St_G(w)=K\}.$$
 \begin{lemma}
 $W_K=\{(a,z)\in W|\Phi_z(a)\in \Irr(L)_K$ is a definable set.
 \end{lemma}
 \begin{proof} Fix a right transversal $(t_i|i=1,\ldots,m)$ for $N$ in $G$ and
 suppose that $K=\cup_{i=1}^s t_iN$.

 Note that $g_i(a)=a^{t_i}$, $a\in L^*$ is a linear
 function of $a$ and $f(a,g)=a^g$, $a\in L^*,g\in L$ is an
 analytic function of $(a,g)$. Define a formula $F_i$ in the language $\La$ in the
 following way.
 $$F_i:=\exists g\ g_i(a)\equiv f(a,g)(\mod z).$$
 Hence $W_K$ is equal to the following definable set
 $$\{(a,z)\in W|F_1(a,z)\&\cdots\& F_s(a,z)\&\neg F_{s+1}(a,z)\&\cdots \& \neg F_m(z,a) \textrm{\ is true in \ } \Z_p
 \}.$$
 \end{proof}
 Now, fix a subgroup $K$, satisfying $N\le K\le G$, and put
 $Q=K/N$. Let $S_Q=\{s_j\}$ be as in the previous section. For
 any
 $w\in \Irr(L)_K$ define $V(w)$ as follows.

 Fix a right transversal $(y_q|q\in Q)$ for $N$ in $K$, such that
$q=y_qN$.
 Since $(a,z)\in K$, we have $(\St_G(w), \St_N(w),\phi) $
 is an extension of $Q$, where $\phi(s)=q$ if $sN=y_qN$ for $s\in
 \St_G(w)$. Put
 $$V(w)=S_Q(\St_G(w), \St_N(w),\phi,w_{\St_N(w)} ).$$
  From Corollary \ref{crit} we know that $V(w)$ takes a
 finite number of values. Let $v$ be a such value. We put
 $$\Irr(L)_{K,v}=\{w\in Irr(L)_K|V(w)=v\}.$$
 \begin{lemma}\label{wkv}
 The set $W_{K,v}=\{(a,z)\in W_K|\Phi_z(a))\in \Irr(L)_{K,v}\}$ is definable.
 \end{lemma}
\begin{proof}
Note that $a(s_i(y_qa_q|q\in Q))$ is an analytic map from
$L^*\times L^{|Q|}$ to $\Z_p$ and  for any $q\in Q$, the function
$g_q(a)=a^{t_q}$, $a\in L^*$ is a linear
 function of $a$. Also the function $f(a,g)=a^g$, $a\in L^*,g\in L$ is an
 analytic function of $(a,g)$.
 Define the
following formula $G_{K,v}$ in $\La$:
$$\begin{array}{lll} G_{K,v}:& = & \exists (n_q|q\in Q)\ \forall
{q\in Q}\ g_q(a)\equiv f(a,-n_q)(\mod
z)\ \& \\
& & ( \forall j \ (a(s_j(y_qn_q))\equiv 0 (\mod z)\ \&\  v=1)\
\bigvee\  \\
& & (a(s_j(y_qn_q))v\equiv 0 (\mod z)\ \& \ a(s_j(y_qn_q))v\not
\equiv 0 (\mod pz))).\end{array}$$ In the first row of the formula
we find a transversal $(y_qn_q|q\in Q)$ for $\St_N(\Phi_z(a))$ in
$\St_G(\Phi_z(a))$ and in the second and third we check the
condition $V(a,z)=v$. Hence we have $W_{K,v}=\{(a,z)\in
W_K|G_{K,v}(a,z)\textrm{\ is true in \ } \Z_p
 \}$. Hence $W_{K,v}$ is definable.
\end{proof}
Now, we are ready to prove the following theorem.

\begin{theorem}
Let $p>2$ and $G$ be a $p$-adic analytic perfect pro-$p$ group.
Then $\zeta^{ch}_G(s)$ is a rational function in $p^s$.
\end{theorem}
\begin{proof}
We conserve the previous notation. For any $N\le K\le G$ we define
$$\zeta^{ch}_{G,K}(s)=\sum_{\lambda \in \Irr(G|\chi),\ \chi\in \Irr(N),\  I_G(\chi)=K}
|\lambda(1)|^{-s}.$$ Note that if $\lambda \in
\Irr(G|\chi),\chi\in \Irr(N)$ and $I_G(\chi)=K$ then $\lambda$
lies over $|G/K|$ irreducible characters of $N$. Hence we have
$$\zeta^{ch}_G(s)=\sum_{N\le K\le
G}\frac{1}{|G:K|}\zeta^{ch}_{G,K}(s).$$ So, we should prove the
rationality of $\zeta^{ch}_{G,K}(s)$ for each $K$.
 Note that
\begin{equation}\label{eq1}\zeta^{ch}_{G,K}(s)=|G/K|^{-s}\sum _{\chi\in \Irr(N),\
I_G(\chi)=K}f_\chi |\chi(1)|^{-s}, \end{equation} where
$f_\chi=\sum_{\lambda \in \Irr(K|\chi)}\left |\frac{\lambda(1)}
{\chi(1)}\right|^{-s}$.

Consider an arbitrary character $\chi \in \Irr(N)$ such that
$I_G(\chi)=K$. Then, from Theorem \ref{charunif}, it follows that
there exists
  a $G$-orbit $\Omega$  in $\Irr(L)$, such that
$\chi=\Phi_{\Omega}$. Let $w\in \Omega$. Since $I_G(\chi)=K$, we
have $w\in \Irr(L)_K$. By Theorem \ref{triplesreduction},
$(K,N,\chi)$ is isomorphic to $(\St_G(w),\St_N(w),w_{\St_N(w)})$.
Hence
$$f_\chi=\sum_{\lambda \in
\Irr(\St_G(w)|w_{\St_N(w)})}{|\lambda(1)}|^{-s}.$$ For $w\in \Irr
(L)$ define $f_w=\sum_{\lambda \in
\Irr(\St_G(w)|w_{\St_N(w)})}{|\lambda(1)|}^{-s}$. Then the
equality (\ref{eq1}) can be rewritten as
$$\zeta^{ch}_{G,K}(s)=|G/K|^{-s}\sum _{w\in \Irr(L)_K,}f_w
|L:\Rad(w)|^{-(s-2)/2}.$$ Now, note that if $w_1, w_2\in
\Irr(L)_{K,v}$ then, by Corollary \ref{crit}, $(\St_G(w_1),
\St_N(w_1),(w_1)_{\St_N(w_1)})$ and  $(\St_G(w_2),
\St_N(w_2),(w_2)_{\St_N(w_2)})$ are isomorphic and, in particular,
$f_{w_1}=f_{w_2}$. Hence in order to prove the rationality of
$\zeta^{ch}_{G,K}(s)$ it is enough to prove that $\sum _{w\in
\Irr(L)_{K,v}} |L:\Rad(w)|^{-s}$ is a rational function in
$p^{-2s}$.

We do it in the same way as we proved Theorem \ref{unif}, because
we have that
$$\sum _{w\in \Irr(L)_{K,v}}
|L:\Rad(w)|^{-s}=p(p-1)^{-1}\int_{W_{K,v}}
 |z|^{-(n+1)}|\alpha(a,z)|^{s}dadz.$$
 The last integral we can transform to a definable integral using the argument of the proof of Theorem \ref{unif},
  because $W_{K,v}$ is a definable set by Lemma \ref{wkv}.
  \end{proof}

\end{document}